\documentclass[12pt,thmsa,sw20bams]{article}

\begin{document}

\begin{titlepage}
\title{\bf  A note over para complex torsion-free affine connection on tangent bundle}
\author{ Mehmet Tekkoyun \footnote{tekkoyun@pau.edu.tr} \\
 {\small Department of Mathematics, Pamukkale University,}\\
{\small 20070 Denizli, Turkey}\\
Ali G\"{o}rg\"{u}l\"{u} \footnote{agorgulu@ogu.edu.tr} \\
 {\small Department of Mathematics, Eski\c{s}ehir Osmangazi University,}\\
 {\small 26480  Eski\c{s}ehir, Turkey}}
\maketitle

\begin{abstract}

The goal of this paper is to introduce the lifting theory that has
an important role in geometry. Therefore, using the lifts of
differential geometric structures we show that tangent bundle $TM$
of paracomplex manifold $M$ admits para-complex torsion-free
affine connection.


\end{abstract}

\section{Introduction}

In differential geometry, lifting method has an important tool. So, using
lift function it may be possible to generalize to differentiable structures
on any manifold or space to its extensions. There are many books and studies
about lift theory. Some of them are given in [1-7]. Lifts of differential
geometric elements defined on any manifold $M$ to tangent manifold $TM$ has
been obtained by Yano and Ishihara [7]. Para-Complex geometry are introduced
by Sch\"{a}fer [8] and, Cruceanu and others [9]. Complex and paracomplex
lift analogues of the geometric structures had been introduced by Tekkoyun
[3,4,5] and Civelek [3,4]. Also, complex and paracomplex lift analogues of
the Lagrangian and Hamiltonian systems in classical mechanics were made by
Tekkoyun and G\"{o}rg\"{u}l\"{u} [6]. In this study, firstly, it is recall
vertical and complete lifts of fundamental structures in geometry. Then, we
deduce that $TM$, tangent bundle of paracomplex manifold $M,$ admits
para-complex torsion-free affine connection.

Throughout this paper, all maps will be understood to be differentiable of
class $C^{\infty }$ and the sum is taken over repeated indices. Also, the
indices $\alpha ,\beta $ are assumed $1\leq \alpha ,\beta \leq m$.

\subsection{Paracomplex Geometry}

A tensor field $J$ of type (1,1) on $M$ such that $J^{2}=I$ is called an
\textit{almost product structure} $J$ on 2m-dimensional manifold $M$. Then,
\textit{almost product manifold }is said to be\textit{\ }the pair $(M,J).$
An \textit{almost paracomplex manifold} is an almost product manifold $(M,J)$
such that the two eigenbundles $T^{\pm }M$ associated to the eigenvalues $%
\pm $1 of $J$, respectively, have the same rank. The dimension of an almost
paracomplex manifold is necessarily even. Equivalently, a splitting of the
tangent bundle $TM$ of manifold $M$, into the Whitney sum of two subbundles
on $T^{\pm }M$ of the same fiber dimension is called an \textit{almost
paracomplex structure} on $M.$ An almost paracomplex structure on manifold $%
M $ may alternatively be defined as a $G$- structure on $M$ with structural
group $GL(n,\mathbf{R})\times GL(n,\mathbf{R})$.

If the $G$- structure defined by the tensor field $J$ is integrable, we call
that an almost paracomplex manifold $(M,J)$ is a \textit{paracomplex
manifold.}

Let $(\ x^{\alpha },\,\ y^{\alpha })$ be a real coordinate system on a
neighborhood $U$ of any point $p$ of $M,$ and $\{(\frac{\partial }{\partial
x^{\alpha }})_{p},(\frac{\partial }{\partial y^{\alpha }})_{p}\}$ and $%
\{(dx^{\alpha })_{p},(dy^{\alpha })_{p}\}$ natural bases over $\mathbf{R}$
of the tangent space and the cotangent space $T_{p}M$ and $T_{p}^{*}M$ of $%
M, $ respectively$.$ Then we explain as
\[
J(\frac{\partial }{\partial x^{\alpha }})=\frac{\partial }{\partial
y^{\alpha }},\,J(\frac{\partial }{\partial y^{\alpha }})=\frac{\partial }{%
\partial x^{\alpha }}
\]
and
\[
J^{*}(dx^{\alpha })=-dy^{\alpha },\,J^{*}(dy^{\alpha })=-dx^{\alpha }.
\]
Let $\ z^{\alpha }=x^{\alpha }+$\textbf{j}$y^{\alpha },\,\overline{z}%
^{\alpha }=x^{\alpha }-$\textbf{j}$y^{\alpha },\,\,$\textbf{j}$^{2}=1,$ be a
paracomplex local coordinate system on a neighborhood $U$ of any point $p$
of $M.$ We express the vector fields and dual covector fields as:
\[
(\frac{\partial }{\partial z^{\alpha }})_{p}=\frac{1}{2}\{(\frac{\partial }{%
\partial x^{\alpha }})_{p}-\mathbf{j}(\frac{\partial }{\partial y^{\alpha }}%
)_{p}\},\,(\frac{\partial }{\partial \overline{z}^{\alpha }})_{p}=\frac{1}{2}%
\{(\frac{\partial }{\partial x^{\alpha }})_{p}+\mathbf{j}(\frac{\partial }{%
\partial y^{\alpha }})_{p}\},
\]
\[
\left( dz^{\alpha }\right) _{p}=\left( dx^{\alpha }\right) _{p}+\mathbf{j}%
(dy^{\alpha })_{p},\,\left( d\overline{z}^{\alpha }\right) _{p}=\left(
dx^{\alpha }\right) _{p}-\mathbf{j}(dy^{\alpha })_{p}.
\]
which symbolize the bases of the tangent space and cotangent space $T_{p}M$
and $T_{p}^{*}M$ of $M$, respectively. Then, using \textbf{j}$^{2}=1$ it is
found
\[
J(\frac{\partial }{\partial z^{\alpha }})=-\mathbf{j}\frac{\partial }{%
\partial z^{\alpha }},\,J(\frac{\partial }{\partial \overline{z}^{\alpha }})=%
\mathbf{j}\frac{\partial }{\partial \overline{z}^{\alpha }}.
\]
The dual map $J^{*}$ of the cotangent space $T_{p}^{*}M$ of manifold $M$ at
any point $p$ satisfies $J^{*2}=I.$ Thus, using \textbf{j}$^{2}=1,$ it is
computed by
\[
J^{*}(dz^{\alpha })=-\mathbf{j}dz^{\alpha },\,J^{*}(d\overline{z}^{\alpha })=%
\mathbf{j}d\overline{z}^{\alpha }.
\]
If vector space $T_{p}M$ is the\thinspace set of tangent vectors $%
Z_{p}=Z^{\alpha }(\frac{\partial }{\partial z^{\alpha }})_{p}+\overline{Z}%
^{\alpha }(\frac{\partial }{\partial \overline{z}^{\alpha }})_{p},$for each $%
\,p\in M,$ then $TM\,$ is the union of $T_{p}M$. Thus, tangent bundle of a
paracomplex manifold $M$ is $\left( TM,\tau _{M},M\right) ,$ where canonical
projection $\tau _{M}$ is $\tau _{M}:TM\rightarrow M$ $(\,\tau
_{M}(\,Z_{p})=p)$ and, also this map are surjective submersion. After now,
coordinates $\left\{ z_{{}}^{\alpha },\overline{z}_{{}}^{\alpha },z_{{}}^{%
\acute{}\alpha },\overline{z}_{{}}^{\acute{}\alpha }\right\} $ are local
coordinates for $TM$.

\section{Lifting theory of paracomplex structures}

\subsection{Lifts of function}

The function $f^{v}\in \mathcal{F}(TM)$ given by
\[
f^{v}=f\circ \tau _{M}
\]

is called \textit{vertical\thinspace \thinspace lift} of\thinspace
\thinspace paracomplex function $f$ $\in \mathcal{F}(M)$\thinspace to $TM$,
where $\tau _{M}:TM\rightarrow M$ canonical projection.$\,\,$\thinspace We
get $rang(f^{v})=rang(f),$ since$\,\,$
\[
f^{v}(Z_{p})=f(\tau _{M}(Z_{p}))=f(p),\,\,\,\,\,\,\,\forall \,Z_{p}\in TM.
\]

The \textit{complete \thinspace lift} of\thinspace \thinspace \thinspace
paracomplex function $f$ \thinspace $\in \mathcal{F}(M)$\ to $TM$ is the
function $f^{c}\in \mathcal{F}(TM)$ given by
\[
f^{c}=z^{\acute{}\alpha }(\frac{\partial f}{\partial z^{\alpha }})^{v}+%
\overline{z}^{\acute{}\alpha }(\frac{\partial f}{\partial \overline{z}%
^{\alpha }})^{v},
\]
where $(z_{{}}^{\alpha },\overline{z}_{{}}^{\alpha },z_{{}}^{\acute{}\alpha
},\overline{z}_{{}}^{\acute{}\alpha })$ are the local coordinates of a
chart-domain $TU\subset TM.$ Further, for $\,Z_{p}\in TM$ it holds
\[
f^{c}(Z_{p})=z^{\acute{}\alpha }(Z_{p})(\frac{\partial f}{\partial z^{\alpha
}})^{v}(p)+\overline{z}^{\acute{}\alpha }(Z_{p})(\frac{\partial f}{\partial
\overline{z}^{\alpha }})^{v}(p).
\]

The general properties of the vertical and complete lifts of paracomplex
functions are as follows:

\[
\begin{array}{ll}
i) & (f.g)^{v}=f^{v}.g^{v},(f+g)^{v}=f^{v}+g^{v}, \\
ii) & (f.g)^{c}=f^{c}.g^{v}+f^{v}.g^{c},(f+g)^{c}=f^{c}+g^{c},
\end{array}
\]

for all $\,f,g\in \mathcal{F} (M).$

\subsection{Lifts of vector field}

In this subsection, we assume that vector field $Z$ is $Z=Z^{\alpha }\frac{%
\partial }{\partial z^{\alpha }}+\overline{Z}^{\alpha }\frac{\partial }{%
\partial \overline{z}^{\alpha }}.$ The vector field $Z^{v}\in \chi (TM)$
determined by
\[
Z^{v}(f^{c})=(Zf)^{v},\,\,\,\,\,\,\forall f\in \mathcal{F}(M)
\]

is the \textit{vertical lift} of a vector field\textbf{\ }$Z\in \chi (M)$ to
$TM.$ Then we get
\[
Z^{v}=(Z^{\alpha })^{v}\frac{\partial }{\partial z^{\acute{}\alpha }}+(%
\overline{Z}^{\alpha })^{v}\frac{\partial }{\partial \overline{z}^{\acute{}%
\alpha }}.
\]

The \textit{complete lift} of a vector field\textbf{\ }$Z\in \chi (M)$ to $%
TM $ is the vector field $Z^{c}\in \chi (TM)$ given by
\[
Z^{c}(f^{c})=(Zf)^{c},{\,\,\,\,\,\,}\forall f\in \mathcal{F}(M).
\]
Clearly, it is gotten
\[
Z^{c}=(Z^{\alpha })^{v}\frac{\partial }{\partial z^{\alpha }}+(\overline{Z}%
^{\alpha })^{v}\frac{\partial }{\partial \overline{z}^{\alpha }}+(Z^{\alpha
})^{c}\frac{\partial }{\partial z^{\acute{}\alpha }}+(\overline{Z}^{\alpha
})^{c}\frac{\partial }{\partial \overline{z}^{\acute{}\alpha }}.
\]

The vertical and complete lifts of paracomplex vector fields have the
following general properties:

\[
\begin{array}{ll}
i) & \,\,(X+Y)^{v}=X^{v}+Y^{v},(X+Y)^{c}=X^{c}+Y^{c}, \\
ii) & (fX)^{v}=f^{v}X^{v},(fX)^{c}=f^{c}X^{v}+f^{v}X^{c}, \\
iii) &
X^{v}(f^{v})=0,X^{c}(f^{v})=X^{v}(f^{c})=(Xf)^{v},X^{c}(f^{c})=(Xf)^{c}, \\
iv) & \left[ X^{v},Y^{v}\right] =0,\left[ X^{v},Y^{c}\right] =\left[
X^{c},Y^{v}\right] =\left[ X,Y\right] ^{v},\left[ X^{c},Y^{c}\right] =\left[
X,Y\right] ^{c} \\
v) & \;\,(\frac{\partial }{\partial z^{\alpha }})^{c}=\frac{\partial }{%
\partial z^{\alpha }},(\frac{\partial }{\partial \overline{z}^{\alpha }}%
)^{c}=\frac{\partial }{\partial \overline{z}^{\alpha }},(\frac{\partial }{%
\partial z^{\alpha }})^{v}=\frac{\partial }{\partial z^{\acute{}\alpha }},(%
\frac{\partial }{\partial \overline{z}^{\alpha }})^{v}=\frac{\partial }{%
\partial \overline{z}^{\acute{}\alpha }},
\end{array}
\]

for all\thinspace \thinspace \thinspace $f\in \mathcal{F}(M),$ $X,Y,Z\in
\chi (M),$ $\,\chi (U)=Sp\left\{ \frac{\partial }{\partial z^{\alpha }},%
\frac{\partial }{\partial \overline{z}^{\alpha }}\right\} ,\,\,\,\chi
(TU)=Sp\left\{ \frac{\partial }{\partial z^{\alpha }},\frac{\partial }{%
\partial \overline{z}^{\alpha }},\frac{\partial }{\partial z^{\acute{}\alpha
}},\frac{\partial }{\partial \overline{z}^{\acute{}\alpha }}\right\} $.

\subsection{Lifts of 1-form}

In this subsection, we consider that the paracomplex 1-form $\omega $ is $%
\omega =\omega _{\alpha }dz^{\alpha }+\overline{\omega }_{\alpha }d\overline{%
z}^{\alpha }.$ The 1-form $\omega ^{v}\in \chi ^{*}(TM)$ explained by
\[
\omega ^{v}(Z^{c})=(\omega Z)^{v},{\,\,\,\,\,\,}\forall Z\in \chi (M),
\]

is said to be the \textit{vertical lift} of a 1-form\textbf{\ }$\omega \in
\chi ^{*}(M)$ to $TM.$

Then, we have
\[
\omega ^{v}=(\omega _{\alpha })^{v}dz^{\alpha }+(\overline{\omega }_{\alpha
})^{v}d\overline{z}^{\alpha }.
\]

The \textit{complete lift} of a 1-form\textbf{\ }$\omega \in \chi ^{*}(M)$
to $TM$ is the 1-form $\omega ^{c}\in \chi ^{*}(TM)$ given by
\[
\omega ^{c}(Z^{c})=(\omega Z)^{c},{\thinspace \thinspace \thinspace
\thinspace \thinspace \thinspace }\forall Z\in \chi (M).
\]

Hence, we compute
\[
\omega ^{c}=(\omega _{\alpha })^{c}dz^{\alpha }+(\overline{\omega }_{\alpha
})^{c}d\overline{z}^{\alpha }+(\omega _{\alpha })^{v}dz^{\acute{}\alpha }+(%
\overline{\omega }_{\alpha })^{v}d\overline{z}^{\acute{}\alpha }.
\]

The properties of the vertical and complete lifts of paracomplex 1-forms are
as follows:

\[
\begin{array}{ll}
i) & (f\omega )^{v}=f^{v}\omega ^{v},(f\omega )^{c}=f^{c}\omega
^{v}+f^{v}\omega ^{c}, \\
ii) & (\omega +\theta )^{v}=\omega ^{v}+\theta ^{v},(\omega +\theta
)^{c}=\omega ^{c}+\theta ^{c} \\
iii) & \omega ^{v}(Z^{v})=0,\omega ^{c}(Z^{v})=\omega ^{v}(Z^{c})=(\omega
Z)^{v},\omega ^{c}(Z^{c})=(\omega Z)^{c}, \\
iv) & (dz^{\alpha })^{c}=\overline{d}z^{\acute{}\alpha },(d\overline{z}%
^{\alpha })^{c}=\overline{d}\overline{z}^{\acute{}\alpha },(dz^{\alpha
})^{v}=\overline{d}z^{\alpha },(d\overline{z}^{\alpha })^{v}=\overline{d}%
\overline{z}^{\alpha },
\end{array}
\]

for all\thinspace \thinspace \thinspace $f\in \mathcal{F}(M)\,,$ $Z\in \chi
(M)$, $\,\omega ,\theta \in \chi ^{*}(M)$, $\,\,\chi ^{*}(U)=Sp\left\{
dz^{\alpha },d\overline{z}^{\alpha }\right\} ,$ $\chi ^{*}(TU)=Sp\left\{
dz^{\alpha },d\overline{z}^{\alpha },dz^{\acute{}\alpha },d\overline{z}^{%
\acute{}\alpha }\right\} ,$ $\overline{d}\,\,$denotes the differential
operator on$\,\,TM$.

\subsection{Lifts of tensor fields of type (1,1)}

The \textit{complete lift} of a paracomplex tensor field of type (1,1)%
\textbf{\ }$F\in \Im _{1}^{1}(M)$ to $TM$ is the tensor field $F^{c}\in \Im
_{1}^{1}(TM)$ given by
\[
F^{c}(Z^{c})=(FZ)^{c},{\,\,\,\,\,\,}\forall Z\in \chi (M).
\]

The complete lift of the paracomplex tensor field of type (1,1)\textbf{\ }$F$
is
\[
\begin{array}{l}
F^{c}=(F_{\alpha }^{\beta })^{v}\frac{\partial }{\partial z^{\beta }}\otimes
dz^{\alpha }+(F_{\alpha }^{\beta })^{c}\frac{\partial }{\partial z^{\acute{}%
\beta }}\otimes dz^{\alpha }+(F_{\alpha }^{\beta })^{v}\frac{\partial }{%
\partial z^{\acute{}\beta }}\otimes dz^{\acute{ }\alpha } \\
\,\,\,\,\,\,\,\,\,\,\,+(\overline{F}_{\alpha }^{\beta })^{v}\frac{\partial }{%
\partial \overline{z}^{\beta }}\otimes d\overline{z}^{\alpha }+( \overline{F}%
_{\alpha }^{\beta })^{c}\frac{\partial }{\partial \overline{z}^{\acute{}%
\beta }}\otimes d\overline{z}^{\alpha }+(\overline{F} _{\alpha }^{\beta
})^{v}\frac{\partial }{\partial \overline{z}^{\acute{}\beta }}\otimes d%
\overline{z}^{\acute{}\alpha }.
\end{array}
\]

\subsection{Lift of para-complex structure}

The \textit{complete lift} of $J$ being a paracomplex tensor field of type
(1,1) is

\[
J^{c}=\mathbf{j}\frac{\partial }{\partial z^{\alpha }}\otimes dz^{\alpha }+%
\mathbf{j}\frac{\partial }{\partial z^{\acute{}\alpha }}\otimes dz^{\acute{}%
\alpha }-\mathbf{j}\frac{\partial }{\partial \overline{z}^{\alpha }}\otimes d%
\overline{z}^{\alpha }-\mathbf{j}\frac{\partial }{\partial \overline{z}^{%
\acute{}\alpha }}\otimes d\overline{z}^{\acute{}\alpha }.
\]

Because of $(J^{c})^{2}=I,$ $J^{c}$ is an almost paracomplex structure for
tangent bundle $TM$.

\section{Para-complex torsion-free affine connection on tangent bundle}

In this section, we assume that $M$ is an almost paracomplex manifold and $%
TM $ its tangent bundle. Let $Z$, $W$ be vector fields and $\nabla $
paracomplex connection, $[\,,]$ Lie bracket on $M$.

The torsion tensor $T$ on $M$ is defined as

\[
T(Z,W)=\nabla _{Z}^{{}}W^{{}}-\nabla _{W^{{}}}^{{}}Z^{{}}-[Z,W^{{}}].
\]

The torsion-free tensor $T^{c}$ on $TM$ determined as
\[
T^{c}(Z^{c},W^{c})=\nabla _{Z^{c}}^{c}W^{c}-\nabla
_{W^{c}}^{c}Z^{c}-[Z^{c},W^{c}]=0
\]

is called the \textit{complete lift} of $T$ on $M.$

The Nijenhuis tensor $N_{J}$ endowed with paracomplex structure $J$ on $M$
is defined as

\[
N_{J}(Z,W)=\left[ Z,W\right] -J\left[ JZ,W\right] -J\left[ Z,JW\right]
+\left[ JZ,JW\right]
\]

Using almost paracomplex structure $J^{c}$ on $TM,$ the \textit{complete lift%
} of $N_{J}$ is the Nijenhuis tensor $N_{J^{c}}^{c}$of $J^{c}$ and given by
\[
N_{J^{c}}^{c}(Z^{c},W^{c})=\left[ Z^{c},W^{c}\right] -J^{c}\left[
J^{c}Z^{c},W^{c}\right] -J^{c}\left[ Z^{c},J^{c}W^{c}\right] +\left[
J^{c}Z^{c},J^{c}W^{c}\right] .
\]

\textbf{Theorem 1: }Let $M$ be almost para-complex manifold and $TM$ its
tangent bundle.\textbf{\ }Every $TM$ fixed with para-complex structure $%
J^{c} $ admits an almost para-complex affine connection with torsion

\[
N_{J^{c}}^{c}=-4T^{c}
\]

where $N_{J^{c}}^{c}$ is the Nijenhuis-tensor of almost para- complex
structure $J^{c}$ and $T^{c}$ is complete lift of the torsion tensor $T$.

\textbf{Proof: }Let\textbf{\ }$\nabla ^{c}$ be torsion-free connection on $%
TM $. We explain $Q^{c}\in \Gamma ((T(T^{*}M))^{2}\otimes T(TM))$ as:

\[
4Q^{c}(X^{c},Y^{c}):=[(\nabla _{J^{c}Y^{c}}^{c}J^{c})X^{c}+J^{c}((\nabla
_{Y^{c}}^{c}J^{c})X^{c})+2J^{c}((\nabla _{X^{c}}^{c }J^{c})Y^{c})]
\]

and furthermore

\[
\widetilde{\nabla }_{X^{c}}^{c}Y^{c}=\nabla _{X^{c}}^{c
}Y^{c}+Q^{c}(X^{c},Y^{c}).
\]

In this case, we calculate

\begin{eqnarray*}
(\widetilde{\nabla }_{X^{c}}^{c}J^{c})Y^{c} &=&\widetilde{\nabla }%
_{X^{c}}^{c}J^{c}Y^{c}-J^{c}\widetilde{\nabla }_{X^{c}}^{c}Y^{c} \\
&=&\nabla _{X^{c}}^{c}J^{c}Y^{c}+Q^{c}(X^{c},J^{c}Y^{c})-J^{c}\nabla
_{X^{c}}^{c}Y^{c}-J^{c}Q^{c}(X^{c},Y^{c}) \\
&=&(\nabla
_{X^{c}}^{c}J^{c})Y^{c}+Q^{c}(X^{c},J^{c}Y^{c})-J^{c}Q^{c}(X^{c},Y^{c})
\end{eqnarray*}

Thinking $A(X^{c},Y^{c})=Q^{c}(X^{c},J^{c}Y^{c})-J^{c}Q^{c}(X^{c},Y^{c}),$
we get

\[
(\widetilde{\nabla }_{X^{c}}^{c}J^{c})Y^{c}=(\nabla
_{X^{c}}^{c}J^{c})Y^{c}+A(X^{c},Y^{c}).
\]

After this, we have to show $A(X^{c},Y^{c})=-(\nabla _{X^{c}}^{c}J^{c})Y^{c}$%
. Also, we write

\begin{eqnarray*}
4Q^{c}(X^{c},J^{c}Y^{c}) &:&=[(\nabla _{Y^{c}}^{c}J^{c})X^{c}+J^{c}((\nabla
_{J^{c}Y^{c}}^{c}J^{c})X^{c})+2J^{c}((\nabla _{X^{c}}^{c}J^{c})J^{c}Y^{c})]
\\
4J^{c}Q^{c}(X^{c},Y^{c}) &:&=[(J^{c}\nabla
_{J^{c}Y^{c}}^{c}J^{c})X^{c}+((\nabla _{Y^{c}}^{c}J^{c})X^{c})+2((\nabla
_{X^{c}}^{c}J^{c})Y^{c})].
\end{eqnarray*}

Taking care of $(J^{c})^{2}=I$ and using $J^{c}((\nabla
_{X^{c}}^{c}J^{c})J^{c}Y^{c})=-(\nabla _{X^{c}}^{c}J^{c})Y^{c}$ and at last
we compute

\[
4A(X^{c},Y^{c})=4Q^{c}(X^{c},J^{c}Y^{c})-4J^{c}Q^{c}(X^{c},Y^{c})=-4(\nabla
_{X^{c}}^{c}J^{c})Y^{c}.
\]

Now, we find the torsion of $\widetilde{\nabla }^{c}$ given as follow:

\[
T^{c^{\widetilde{\nabla }^{c}}}(X^{c},Y^{c})=T^{c^{\nabla
^{c}}}(X^{c},Y^{c})+Q^{c}(X^{c},Y^{c})-Q^{c}(Y^{c},X^{c})=Q^{c}(X^{c},Y^{c})-Q^{c}(Y^{c},X^{c}).
\]

Considering the definition of $Q^{c}$ we have

\begin{eqnarray*}
4T^{c^{\widetilde{\nabla }^{c}}}(X^{c},Y^{c})
&=&4Q^{c}(X^{c},Y^{c})-4Q^{c}(Y^{c},X^{c}) \\
&=&[(\nabla _{J^{c}Y^{c}}^{c}J^{c})X^{c}+J^{c}((\nabla
_{Y^{c}}^{c}J^{c})X^{c})+2J^{c}((\nabla _{X^{c}}^{c}J^{c})Y^{c})] \\
&&-[(\nabla _{J^{c}X^{c}}^{c}J^{c})Y^{c}+J^{c}((\nabla
_{X^{c}}^{c}J^{c})Y^{c})+2J^{c}((\nabla _{Y^{c}}^{c}J^{c})X^{c})]
\end{eqnarray*}

Making necessary operations and taking $[X^{c},Y^{c}]=\nabla
_{X^{c}}^{c}Y^{c}-\nabla _{Y^{c}}^{c}X^{c},(J^{c})^{2}=I,$ finally it is
shown

\begin{eqnarray*}
4T^{c^{\widetilde{\nabla }^{c}}}(X^{c},Y^{c})
&=&-[J^{c}X^{c},J^{c}Y^{c}]-[X^{c},Y^{c}]+J^{c}[X^{c},J^{c}Y^{c}]+J^{c}[X^{c},J^{c}Y^{c}]
\\
&=&-N_{J^{c}}^{c}(X^{c},Y^{c}).
\end{eqnarray*}

\textbf{Corollary: }Every tangent bundle $TM$ endowed with para-complex
structure $J^{c}$ admits a para-complex torsion-free affine connection.

\thinspace \textbf{REFERENCES}\thinspace

\thinspace \thinspace \thinspace \thinspace [1] De Leon M., Rodrigues P.R.,
Methods of Differential Geometry in Analytical Mechanics, North-Hol. Math.
St.,152, Elsevier Sc. Pub. Com., Inc.,Amsterdam, 1989.

[2] Etoya, J.J., On a Complete Lifting of Derivations, Tensor, \textbf{38},
(1982)169-178.

[3] Tekkoyun M., Civelek \c{S}., First Order Lifts of Complex Structures,
Algebras Groups and Geometries (AGG), \textbf{19}(3), (2002)373-382.

[4] Tekkoyun M., Civelek \c{S}., On Lifts of Structures on Complex
Manifolds, Differential Geometry-Dynamics Systems, \textbf{5}(1),
(2003)59-64.

[5] Tekkoyun M., On Lifts of Paracomplex Structures, Turk. J. Math., \textbf{%
30}, (2006)197-210.

[6] Tekkoyun M., G\"{o}rg\"{u}l\"{u} A., Higher Order Complex Lagrangian and
Hamiltonian Mechanics Systems , Physics Letters A, \textbf{357(}2006)261-269.

[7] Yano, K., Ishihara, S., Tangent and Cotangent Bundles, Marcel Dekker
Inc., New York, 1973.

[8] Sch\"{a}fer L., tt*-Bundles in Para-Complex Geometry, Special Para-K\"{a}%
hler Manifolds and Para-Pluriharmonics Maps, Supported by a grant of the
'Studientstiftung des deutschen Volkes'.

[9] V. Cruceanu, P.M. Gadea, J. M. Masqu\'{e}, Para-Hermitian and Para-
K\"{a}hler Manifolds, Supported by the commission of the European
Communities Action for Cooperation in Sciences and Technology with Central
Eastern European Countries n. ERB3510PL920841
\end{titlepage}

\end{document}